\documentclass[12pt]{amsart}

\usepackage{amsmath}
\usepackage{amssymb}

\theoremstyle{definition}

\theoremstyle{remark}

\theoremstyle{plain}

\numberwithin{equation}{section}

\def\CC{{\mathbb C}}

\def\QQ{{\mathbb Q}}
\def\RR{{\mathbb R}}
\def\TT{{\mathbb T}}
\def\ZZ{{\mathbb Z}}

\def\H{{\mathfrak H}}

\def\e{\mathrm{e}}
\def\i{\mathrm{i}}

\def\id{\operatorname{id}}

\def\C{\operatorname{C{}}}

\def\L{\operatorname{L{}}}

\def\Sw{{\mathcal S}}
\def\SL{\operatorname{SL}}

\def\vecl{{\text{\boldmath$l$}}}
\def\vecm{{\text{\boldmath$m$}}}
\def\vecn{{\text{\boldmath$n$}}}

\def\vecs{{\text{\boldmath$s$}}}

\def\vecw{{\text{\boldmath$w$}}}
\def\vecx{{\text{\boldmath$x$}}}
\def\vecy{{\text{\boldmath$y$}}}
\def\vecalf{{\text{\boldmath$\alpha$}}}

\def\vecxi{{\text{\boldmath$\xi$}}}

\def\vecnull{{\text{\boldmath$0$}}}

\def\scrs{{\mathcal S}}

\def\Im{\operatorname{Im}}

\begin{document}

\title[Inhomogeneous quadratic forms II]
{Pair correlation densities of inhomogeneous quadratic forms II}
\author{Jens Marklof}
\address{School of Mathematics, University of Bristol,
Bristol BS8 1TW, U.K.} 
\address{{\tt j.marklof@bris.ac.uk}}
\thanks{December 2000/September 2002, to appear in Duke Math.~Journal}
\thanks{MSC2000: 11P21 (11F27, 22E40, 58F11)}

\begin{abstract}
Denote by $\| \,\cdot\,\|$ the euclidean norm in $\RR^k$.
We prove that the local pair correlation density
of the sequence $\| \vecm -\vecalf \|^k$, $\vecm\in\ZZ^k$,
is that of a Poisson process, under diophantine conditions
on the fixed vector $\vecalf\in\RR^k$: in dimension two,
vectors $\vecalf$ of any diophantine type are admissible; 
in higher dimensions ($k>2$), Poisson statistics are only observed 
for diophantine vectors of type $\kappa<(k-1)/(k-2)$. 
Our findings support a 
conjecture of Berry and Tabor on the Poisson nature of
spectral correlations in quantized integrable systems.
\end{abstract} 

\maketitle

\section{Introduction}

\subsection{\label{intro}}
Berry and Tabor \cite{Berry77} have conjectured that the
local correlations of quantum energy levels of integrable
systems are those of independent random numbers from a
Poisson process. We will here present a proof of this conjecture
for the two-point correlations of the sequence
$$
0\leq\lambda_1\leq\lambda_2\leq \cdots\rightarrow\infty
$$
given by the values of 
$$
\|\vecm-\vecalf\|^2 = (m_1-\alpha_1)^2 + \cdots + (m_k-\alpha_k)^2
$$
at lattice points $\vecm=(m_1,\ldots,m_k)\in\ZZ^k$, 
for fixed $\vecalf=(\alpha_1,\ldots,\alpha_k)\in\RR^k$.
These numbers represent the eigenvalues of the Laplacian 
$$
-\Delta = -\frac{\partial^2}{\partial x_1^2}
-\ldots -\frac{\partial^2}{\partial x_k^2}
$$
on the flat torus $\TT^k$ with quasi-periodicity conditions
$$
\varphi(\vecx+\vecl) = \e^{-2\pi\i \vecalf\cdot\vecl} \varphi(\vecx) , \quad
\vecl\in\ZZ^k,
$$
and may therefore be viewed as energy levels of the
quantized geodesic flow. Statistical properties of the 
above sequence were first studied by
Cheng, Lebowitz and Major \cite{Cheng91,Cheng94}
in dimension $k=2$.
We will here extend our studies \cite{Marklof00,Marklof00d}
to dimensions $k\geq 2$.

Previous results on the Berry-Tabor conjecture for flat tori include
\cite{Eskin98b,Marklof98,Sarnak97} in dimension
$k=2$ and \cite{VanderKam96,VanderKam97,VanderKam00} 
for $k>2$. For more details and references
see \cite{Bleher99,Marklof98,Marklof00c,Sarnak97I}. 

\subsection{\label{start}}
We are interested in the local correlations between the
$\lambda_j$ on the scale of the mean spacing.
Because the mean density is increasing as $\lambda\rightarrow\infty$, i.e.,
$$
\frac1\lambda \#\{ j : \lambda_j\leq\lambda \}
=\frac1\lambda \#\{ \vecm\in\ZZ^k: \|\vecm-\vecalf\|^2 \leq \lambda \}
\sim B_k \lambda^{k/2-1} ,
$$
where $B_k$ is the volume of the unit ball,
it is necessary to rescale the sequence by setting
$$
X_j = \lambda_j^{k/2} .
$$
Then
$$
\frac1X \#\{ j : X_j\leq X \}
=\frac1X \#\{ \vecm\in\ZZ^k: \|\vecm-\vecalf\|^k \leq X \}
\rightarrow B_k 
$$
for $X \rightarrow\infty$, 
and hence the mean spacing is constant, as required.

\subsection{}
The {\em pair correlation density} of a sequence with 
constant mean density $D$ is defined as
$$
R_2[a,b](X) = \frac{1}{DX}  \#\{ i\neq j : X_i,X_j\in[X,2X],\;
X_i-X_j\in[a,b] \} .
$$

We recall the following classical result.

\subsection{Theorem}{\em
If the $X_j$ come from a Poisson process with mean density $D$,
one has
$$
\lim_{X\rightarrow\infty} R_2[a,b](X)  = D (b-a) 
$$
almost surely.
}

\subsection{}
We will here prove a similar result for the {\em deterministic}
sequence in \ref{intro}, which holds, however, only 
under diophantine conditions on $\vecalf$. 
The vector $\vecalf=(\alpha_1,\ldots,\alpha_k)\in\RR^k$ 
is said to be {\em diophantine of type $\kappa$}, if there exists a 
constant $C$ such that 
$$
\max_j |\alpha_j - \frac{m_j}{q}| > \frac{C}{q^\kappa}
$$
for all $m_1,\ldots,m_k,q\in\ZZ$, $q>0$. The smallest possible value
for $\kappa$ is $\kappa=1+\frac1k$. In this case $\vecalf$
is called {\em badly approximable}.

\subsection{Theorem\label{pairthm}}{\em
Suppose $\vecalf$ is diophantine of type $\kappa<\frac{k-1}{k-2}$
and the components of the vector $(\vecalf,1)\in\RR^{k+1}$
are linearly independent over $\QQ$. 
Then
$$
\lim_{X\rightarrow\infty} R_2[a,b](X)  = B_k (b-a) .
$$
}

The condition in the theorem is satisfied, if for instance
the components of $(\vecalf,1)$ form a basis of a real algebraic
number field of the degree $k+1$. In this case $\kappa=1+\frac1k$
\cite{Schmidt80}.

The condition $\kappa<\frac{k-1}{k-2}$ in Theorem \ref{pairthm}
is sharp: 

\subsection{Theorem\label{sharp}}{\em
Let $k>2$.
For any $a>0$, there exists a set $C\subset\TT^k$ of second Baire category,
for which the following holds.

{\rm (i)} 
All $\vecalf\in C$ are diophantine of type $\kappa=\frac{k-1}{k-2}$,
and the components of the vector $(\vecalf,1)\in\RR^{k+1}$
are linearly independent over $\QQ$. 

{\rm (ii)} For $\vecalf\in C$, we find arbitrarily large $X$
such that
$$
R_2[-a,a](X)  \geq \frac{\log X}{\log\log\log X}.
$$

{\rm (iii)} For $\vecalf\in C$, there exists an infinite  sequence
$L_1<L_2<\cdots\rightarrow\infty$ such that
$$
\lim_{j\rightarrow\infty} R_2[-a,a](L_j)  = 2\pi a .
$$
}

In Theorem \ref{sharp} (ii), $\log\log\log X$ may be replaced by any
slowly increasing positive function $\nu(X)\leq \log\log\log X$
with $\nu(X)\rightarrow\infty$ as $X\rightarrow\infty$.

Without imposing any diophantine condition, the rate of divergence may
be even worse:

\subsection{Theorem\label{except}}{\em
For any $a>0$, there exists a set $C\subset\TT^k$ of second Baire category,
for which the following holds.

{\rm (i)} 
For $\vecalf\in C$, the components of the vector $(\vecalf,1)\in\RR^{k+1}$
are linearly independent over $\QQ$. 

{\rm (ii)} For $\vecalf\in C$, we find arbitrarily large $X$
such that
$$
R_2[-a,a](X)  \geq 
\begin{cases}
\dfrac{\log X}{\log\log\log X} & \qquad (k=2)\\[10pt]
\dfrac{X^{(k-2)/k}}{\log\log\log X} & \qquad (k>2).
\end{cases}
$$

{\rm (iii)} For $\vecalf\in C$, there exists an infinite  sequence
$L_1<L_2<\cdots\rightarrow\infty$ such that
$$
\lim_{j\rightarrow\infty} R_2[-a,a](L_j)  = 2\pi a .
$$
}

Again, $\log\log\log X$ may be replaced by any
slowly increasing positive function $\nu(X)\leq \log\log\log X$
with $\nu(X)\rightarrow\infty$ as $X\rightarrow\infty$.

Theorems \ref{sharp} and \ref{except} are proved in Section \ref{seccounter}.

\section{Rescaling}

\subsection{}
We shall see in this section, how
Theorem \ref{pairthm}, which is the central result of this paper, 
follows as a straightforward corollary from the asymptotics of
the generalized pair correlation function
$$
R_2(\psi,\lambda)
=\frac{1}{B_k \lambda^{k/2}} \sum_{i,j=1}^\infty 
\psi\big(\frac{\lambda_i}{\lambda},\frac{\lambda_j}{\lambda},
\lambda^{k/2-1}(\lambda_i-\lambda_j)\big) ,
$$
with $\psi\in\C_0(\RR^+\times\RR^+\times\RR)$,
i.e., continuous and of compact support.

\subsection{Theorem\label{glthm}}{\em 
Let $\psi\in\C_0(\RR^+\times\RR^+\times\RR)$.
Suppose the components of $(\vecalf,1)\in\RR^{k+1}$ 
are linearly independent over $\QQ$, and assume 
$\vecalf$ is diophantine of type $\kappa<\frac{k-1}{k-2}$.
Then
$$
\lim_{\lambda\rightarrow\infty}
R_2(\psi,\lambda) 
= \frac{k}{2} \int_0^\infty \psi(r,r,0) r^{k/2-1} dr  
+ \frac{k^2}{4} B_k \int_{\RR}\int_0^\infty \psi(r,r,s) r^{k-2} dr \, ds  .
$$
}

\subsection{Theorem \ref{glthm} $\mathbf\Rightarrow$ Theorem \ref{pairthm}.}
Let us now show how Theorem \ref{glthm} implies
Theorem \ref{pairthm}. For $\psi_1,\psi_2\in\C_0(\RR^+)$
with support in the compact interval $I$ not containing the origin 0,
and $\sigma\in\C_0(\RR)$, we define
$$
\psi(r_1,r_2,s) = \psi_1(r_1^{k/2})\psi_2(r_2^{k/2})\sigma(\rho(r_1,r_2) s),
$$
with
$$
\rho(r_1,r_2)=\frac{r_1^{k/2}-r_2^{k/2}}{r_1-r_2}
=
\begin{cases}
\displaystyle
\sum_{\nu=1}^{k/2} r_1^{k/2-\nu} r_2^{\nu-1} & \text{($k$ even)} \\
\displaystyle
\frac{1}{r_1^{1/2}+r_2^{1/2}}
\sum_{\nu=1}^k r_1^{(k-\nu)/2} r_2^{(\nu-1)/2} & \text{($k$ odd)} .
\end{cases}
$$
It is evident that we can find a constant $\delta>0$ such that
$$
\delta<\rho(r_1,r_2)<\frac1\delta
$$
uniformly for all $r_1,r_2\in I$.

The assumptions on $\psi$ in Theorem \ref{glthm} are therefore satisfied, 
giving
\begin{multline*}
\lim_{\lambda\rightarrow\infty}
\frac{1}{B_k \lambda^{k/2}} \sum_{i,j=1}^\infty 
\psi_1(\frac{\lambda_i^{k/2}}{\lambda^{k/2}})
\psi_2(\frac{\lambda_j^{k/2}}{\lambda^{k/2}})
\sigma(\lambda_i^{k/2}-\lambda_j^{k/2}) \\
= \frac{k}{2} \sigma(0)
\int_0^\infty \psi_1(r^{k/2})\psi_2(r^{k/2}) r^{k/2-1} dr  \\
+ \frac{k^2}{4} B_k \int_{\RR}\int_0^\infty \psi_1(r^{k/2})\psi_2(r^{k/2})
\sigma(\rho(r,r)s) r^{k-2} dr \, ds  .
\end{multline*}
With $\rho(r,r)=\tfrac{k}{2} r^{k/2-1}$ and the substitutions  
$X=\lambda^{k/2}$, $x=r^{k/2}$ and $s\mapsto s/\rho(r,r)$ we finally have
\begin{multline*}
\lim_{X\rightarrow\infty}
\frac{1}{B_k X} \sum_{i,j=1}^\infty 
\psi_1(\frac{X_i}{X})
\psi_2(\frac{X_j}{X})
\sigma(X_i-X_j) \\
= \sigma(0) \int_0^\infty \psi_1(x)\psi_2(x)  dx  
+ B_k \int_{\RR} \sigma(s)\, ds \int_0^\infty \psi_1(x)\psi_2(x)\, dx.
\end{multline*}
The first term on the right-hand side comes obviously from the diagonal
terms $X_i=X_j$ (use the asymptotics in \ref{start}), so
$$
\lim_{X\rightarrow\infty}
\frac{1}{B_k X} \sum_{i\neq j} 
\psi_1(\frac{X_i}{X})
\psi_2(\frac{X_j}{X})
\sigma(X_i-X_j) \\
= B_k \int_{\RR} \sigma(s)\, ds \int_0^\infty \psi_1(x)\psi_2(x)\, dx,
$$
which is a smoothed version of Theorem \ref{pairthm}.
We complete the proof by quoting a standard density argument
(compare proof of Theorem 1.8 in \cite{Marklof00}), in which
the characteristic functions of the intervals $[1,2]$, $[1,2]$
and $[a,b]$ are approximated 
from above and below by smooth functions $\psi_1$, $\psi_2$
and $\sigma$, respectively.
\qed

\subsection{}
It will be sufficient to restrict our attention to the following special case
of Theorem \ref{glthm}. Put
$$
R_2(\psi_1,\psi_2,h,\lambda)
=\frac{1}{B_k \lambda^{k/2}} \sum_{i,j=1}^\infty 
\psi_1(\frac{\lambda_i}{\lambda})\psi_2(\frac{\lambda_j}{\lambda})
\hat h\big(\lambda^{k/2-1}(\lambda_i-\lambda_j)\big) ,
$$
Here $\psi_1,\psi_2\in\Sw(\RR_+)$ are 
real-valued, and $\Sw(\RR_+)$ denotes the Schwartz
class of infinitely differentiable functions
of the half line $\RR_+$ (including the origin), which, as well as their
derivatives, decrease rapidly at $+\infty$.
$\hat h$ is the Fourier transform of a compactly supported function
$h\in\C_0(\RR)$,
$$
\hat h(s)= \int_\RR h(u) e(\tfrac12 us)\, du,
$$
with the shorthand $e(z):=\e^{2\pi\i z}$.

We will prove the following (Section \ref{secmain}).

\subsection{Theorem\label{lthm}}{\em 
Let $\psi_1,\psi_2\in\Sw(\RR_+)$
and $h\in\C_0(\RR)$.
Suppose the components of $(\vecalf,1)\in\RR^{k+1}$ 
are linearly independent over $\QQ$, and assume 
$\vecalf$ is diophantine of type $\kappa<\frac{k-1}{k-2}$.
Then
\begin{multline*}
\lim_{\lambda\rightarrow\infty}
R_2(\psi_1,\psi_2,h,\lambda)
= \frac{k}{2} \hat h(0) \int_0^\infty \psi_1(r) \psi_2(r) r^{k/2-1} dr  \\
+ \frac{k^2}{4} B_k \int \hat h(s)\, ds\;
\int_0^\infty \psi_1(r)\psi_2(r) r^{k-2} dr.
\end{multline*}
}

\subsection{Theorem \ref{lthm} $\mathbf\Rightarrow$ Theorem \ref{glthm}.}
For any fixed $\epsilon>0$ we find finite linear combinations
(cf.~Section 8.6 in \cite{Marklof00})
$$
\psi^\pm(r_1,r_2,s)=\sum_\nu \psi_{1,\nu}^\pm(r_1)\psi_{2,\nu}^\pm(r_2) 
\hat h_{\nu}^\pm(s)
$$
of functions satisfying the conditions of Theorem \ref{lthm}
such that
$$
\psi^-(r_1,r_2,s)\leq\psi(r_1,r_2,s) \leq\psi^+(r_1,r_2,s)
$$
and
$$
\iint (\psi^+(r,r,s)-\psi^-(r,r,s)) r^{k-2} dr \, ds
<\epsilon .
$$
Theorem \ref{lthm} tells us that
$$
\lim_{\lambda\rightarrow\infty}
\frac{1}{B_k \lambda^{k/2}} \sum_{i\neq j} 
\psi^\pm\big(\frac{\lambda_i}{\lambda},\frac{\lambda_j}{\lambda},
\lambda^{k/2-1}(\lambda_i-\lambda_j)\big)
=\frac{k^2}{4} B_k \iint \psi^\pm(r,r,s) r^{k-2} dr \, ds 
$$
(recall the first term in that theorem comes trivially from
the diagonal terms $i=j$).
This implies
$$
\limsup_{\lambda\rightarrow\infty}
\frac{1}{B_k \lambda^{k/2}} \sum_{i\neq j} 
\psi\big(\frac{\lambda_i}{\lambda},\frac{\lambda_j}{\lambda},
\lambda^{k/2-1}(\lambda_i-\lambda_j)\big) 
\leq \frac{k^2}{4} B_k (\iint \psi(r,r,s) r^{k-2} dr \, ds +
\epsilon)
$$
and
$$
\liminf_{\lambda\rightarrow\infty}
\frac{1}{B_k \lambda^{k/2}} \sum_{i\neq j} 
\psi\big(\frac{\lambda_i}{\lambda},\frac{\lambda_j}{\lambda},
\lambda^{k/2-1}(\lambda_i-\lambda_j)\big)  
\geq \frac{k^2}{4} B_k (\iint \psi(r,r,s) r^{k-2} dr \, ds 
-\epsilon) .
$$
Because these inequalities hold for arbitrarily small $\epsilon>0$,
Theorem \ref{glthm} must be true.
\qed

\section{Outline of the proof of Theorem \ref{lthm}\label{secoutline}}

Using the Fourier transform we may write
\begin{equation*}
\begin{split}
R_2 & (\psi_1,\psi_2,h,\lambda) \\ 
= & \frac{1}{B_k} \int 
\big(\frac{1}{\lambda^{k/4}}
\sum_{j} \psi_1(\frac{\lambda_j}{\lambda}) e(\tfrac12 \lambda_j 
\lambda^{k/2-1}u) \big) 
 \overline{
\big(\frac{1}{\lambda^{k/4}}
\sum_{j} \psi_2(\frac{\lambda_j}{\lambda}) e(\tfrac12 \lambda_j 
\lambda^{k/2-1}u) \big)}
h(u)\, du \\
= & \frac{1}{B_k\lambda^{k/2-1}} \int 
\big(\frac{1}{\lambda^{k/4}}
\sum_{j} \psi_1(\frac{\lambda_j}{\lambda}) e(\tfrac12 \lambda_j u) \big) 
 \overline{
\big(\frac{1}{\lambda^{k/4}}
\sum_{j} \psi_2(\frac{\lambda_j}{\lambda}) e(\tfrac12 \lambda_j u) \big)}
h(\lambda^{-(k/2-1)}u)\, du .
\end{split}
\end{equation*}
The sum
$$
\theta_\psi(u,\lambda) =
\frac{1}{\lambda^{k/4}}
\sum_{j} \psi(\frac{\lambda_j}{\lambda}) e(\tfrac12 \lambda_j u)
$$ 
will be identified as a Jacobi theta sum 
living on a certain noncompact but finite-volume 
manifold $\Sigma$ (Section \ref{sectheta}). 
The integration in
$$
R_2(\psi_1,\psi_2,h,\lambda) 
= \frac{1}{B_k} \lambda^{-(k/2-1)} \int \theta_{\psi_1}(u,\lambda)
\overline{\theta_{\psi_2}(u,\lambda)} h(\lambda^{-(k/2-1)} u)\, du
$$
amounts to averaging along a unipotent orbit on
$\Sigma$, which becomes equidistributed as $\lambda\rightarrow\infty$
(Section \ref{secequi}). 
Diophantine conditions on $\vecalf$
are necessary to secure the convergence of the limit
(Section \ref{secdiophantine}). 

The equidistribution theorem yields then
$$
\frac{1}{\mu(\Sigma)}\int_\Sigma \theta_{\psi_1}
\overline{\theta_{\psi_2}} d\mu \; \int h(u)\, du ,
$$
where $\mu$ is the invariant measure.
The first integral can be calculated quite easily (Section \ref{secmain}), 
and we will see that
$$
\frac{1}{\mu(\Sigma)}\int_\Sigma  \theta_{\psi_1}
\overline{\theta_{\psi_2}} d\mu \; \int h(u)\, du 
= \frac{k}{2} B_k \int \psi_1(r) \psi_2(r)\,r^{k/2-1} dr \; \int h(u)\, du ,
$$
which finally yields
$$
\frac{k}{2} B_k \hat h(0) \int \psi_1(r) \psi_2(r)\, r^{k/2-1}dr ,
$$
compare the first term in Theorem \ref{lthm}.

An additional contribution comes from an arc of the orbit, which vanishes
into the cusp. Even though the length of that arc tends to zero, 
the average over the unbounded theta function gives a non-vanishing 
contribution
$$
\frac{k^2}{2} B_k^2 h(0) \int \psi_1(r) \psi_2(r)\,r^{k-2} dr 
=  \frac{k^2}{4} B_k^2 \int \hat h(u)\,du \int \psi_1(r) \psi_2(r)\, 
r^{k-2}dr ,
$$
which corresponds to the second term in Theorem \ref{lthm}.

\section{Theta sums\label{sectheta}}

\subsection{}\label{iwasawa}
Consider the semi-direct product group
$G^k=\SL(2,\RR)\ltimes\RR^{2k}$ with multiplication law
$$
(M;\vecxi)(M';\vecxi') =(MM';\vecxi+M\vecxi') ,
$$
where $M,M'\in\SL(2,\RR)$ and $\vecxi,\vecxi'\in\RR^{2k}$;
the action of $\SL(2,\RR)$ on $\RR^{2k}$ is defined canonically as
$$
M\vecxi = \begin{pmatrix}
a\vecx +b\vecy \\
c\vecx +d\vecy \\
\end{pmatrix}, \quad 
M=\begin{pmatrix}
a&b\\
c&d
\end{pmatrix},\quad 
\vecxi= \begin{pmatrix}
\vecx\\
\vecy
\end{pmatrix},
$$
where $\vecx,\vecy\in\RR^k$.
A convenient parametrization of $\SL(2,\RR)$ can be obtained
by means of the Iwasawa decomposition 
$$
M = \begin{pmatrix}
1 & u\\
0 & 1
\end{pmatrix}
\begin{pmatrix}
v^{1/2} & 0\\
0 & v^{-1/2}
\end{pmatrix}
\begin{pmatrix}
\cos\phi & -\sin\phi\\
\sin\phi &  \cos\phi
\end{pmatrix},
$$
which is unique for $\tau=u+\i v\in\H$, $\phi\in[0,2\pi)$,
where $\H$ denotes the upper half plane $\H=\{ \tau \in\CC: \Im\tau>0\}$.

\subsection{}
For any Schwartz function $f\in\Sw(\RR^k)$ we define
the {\em Jacobi theta sum} $\Theta_f$ by
$$
\Theta_f(\tau,\phi;\vecxi) = v^{k/4} 
\sum_{\vecm\in\ZZ^k} 
f_\phi( (\vecm-\vecy) v^{1/2})\, e(\tfrac12\|\vecm-\vecy\|^2 u 
+ \vecm\cdot \vecx) ,
$$
where
$$
f_\phi(\vecw)= \int_{\RR^k} G_\phi(\vecw,\vecw')  f(\vecw') \,dw' ,
$$
with the integral kernel
$$
G_\phi(\vecw,\vecw')=e(-k \sigma_\phi/8) |\sin\phi|^{-k/2}
e\left[\frac{\tfrac12(\|\vecw\|^2+\|\vecw'\|^2)
\cos\phi- \vecw\cdot\vecw'}{\sin\phi}\right],
$$
where $\sigma_\phi=2\nu+1$ when $\nu\pi<\phi<(\nu+1)\pi$, $\nu\in\ZZ$.
The operators $U^\phi: f \mapsto f_\phi$ 
are unitary, see \cite{Lion80,Marklof99} for details.
Note in particular $U^0=\id$.

The proofs of the remaining statements in this section are found 
in Section 4 of \cite{Marklof00}.

\subsection{Lemma}\label{schwartz}
{\em Let $f_\phi=U^\phi f$, with $f\in{\mathcal S}(\RR^k)$. Then, 
for any $R>1$, there is a constant $c_R$ such that for all $\vecw\in\RR^k$,
$\phi\in\RR$, we have
$$
|f_\phi(\vecw)| \leq c_R (1+\|\vecw\|)^{-R} .
$$
}

\subsection{}

Let us consider the following discrete subgroup in $G^k$.
$$
\Gamma^k=\left\{ 
( \begin{pmatrix}
a & b \\
c & d
\end{pmatrix}
; 
\begin{pmatrix}
ab \vecs \\
cd \vecs
\end{pmatrix}
+\vecm )
: \; \begin{pmatrix}
a & b \\
c & d
\end{pmatrix} \in\SL(2,\ZZ),\; \vecm\in\ZZ^{2k}
\right\} \subset G^k,
$$
with
$\vecs=(\tfrac12,\tfrac12,\ldots,\tfrac12)\in\RR^k$. 

\subsection{Lemma}\label{generator}
{\em $\Gamma^k$ is generated by the elements
$$
( \begin{pmatrix}
0 & -1 \\
1 & 0
\end{pmatrix}
; \vecnull ), \quad
( \begin{pmatrix}
1 & 1 \\
0 & 1
\end{pmatrix}
; \begin{pmatrix}
\vecs \\
\vecnull
\end{pmatrix}), \quad
( \begin{pmatrix}
1 & 0 \\
0 & 1
\end{pmatrix}
; \vecm), \quad \vecm\in\ZZ^{2k} .
$$
}

\subsection{Proposition}
{\em The left action of the group $\Gamma^k$ on $G^k$ 
is properly discontinuous.
A fundamental domain of $\Gamma^k$ in $G^k$ is given by
$$
{\mathcal F}_{\Gamma^k} = {\mathcal F}_{\SL(2,\ZZ)}
\times \{ \phi\in[0,\pi) \} \times \{\vecxi\in[-\tfrac12,\tfrac12)^{2k} \}.
$$
where ${\mathcal F}_{\SL(2,\ZZ)}$ 
is the fundamental domain in $\H$ of the modular group
$\SL(2,\ZZ)$, given by $\{ \tau\in\H: u\in[-\tfrac12,\tfrac12), |\tau|> 1 \}$.
}

\subsection{Proposition}
{\em For $f,g\in\Sw(\RR^k)$,
$\Theta_f(\tau,\phi;\vecxi)\overline{\Theta_g(\tau,\phi;\vecxi)}$ 
is invariant under the left action
of $\Gamma^k$.}

\subsection{Proposition\label{asymp}}{\em
Let $f,g\in\Sw(\RR^k)$.  For any $R>1$, we have
$$
\Theta_f(\tau,\phi;\begin{pmatrix}\vecx \\ \vecy \end{pmatrix})
\overline{\Theta_g(\tau,\phi;\begin{pmatrix}\vecx \\ \vecy \end{pmatrix})}
= v^{k/2} 
\sum_{\vecm\in\ZZ^k}
f_\phi((\vecm-\vecy) v^{1/2}) \overline{g_\phi((\vecm-\vecy) v^{1/2})}
+O_R(v^{-R})
$$
uniformly for all $(\tau,\phi;\vecxi)\in G^k$ with $v>\frac12$.
In addition
$$
\Theta_f(\tau,\phi;\begin{pmatrix}\vecx \\ \vecy \end{pmatrix})
\overline{\Theta_g(\tau,\phi;\begin{pmatrix}\vecx \\ \vecy \end{pmatrix})}
= v^{k/2} 
f_\phi((\vecn-\vecy) v^{1/2}) \overline{g_\phi((\vecn-\vecy) v^{1/2})}
+O_R(v^{-R}),
$$
uniformly for all $(\tau,\phi;\vecxi)\in G^k$ with $v>\frac12$,
$\vecy\in\vecn+[-\tfrac12,\tfrac12]^k$ and $\vecn\in\ZZ^k$.
}

\subsection{Lemma}\label{three}
{\em
The subgroup
$$
\Gamma_\theta \ltimes \ZZ^{2k},
$$
where 
$$
\Gamma_\theta=\left\{ 
\begin{pmatrix}
a & b \\
c & d
\end{pmatrix} \in\SL(2,\ZZ):\; ab\equiv cd \equiv 0\bmod 2
\right\} 
$$ 
is the theta group, is of index three in $\Gamma^k$.
}

\subsection{Lemma}\label{two}
{\em 
$\Gamma^k$ is of finite index in $\SL(2,\ZZ)\ltimes (\frac12\ZZ)^{2k}$.
}

\subsection{}
Note: The theta sum defined in this section 
is related to the sum $\theta_{\psi_1}(u,\lambda)$ 
in Section \ref{secoutline} by
$$
\theta_{\psi_1}(u,\lambda)\overline{\theta_{\psi_2}(u,\lambda)}
=\Theta_f(u+\i\frac1\lambda,0;\begin{pmatrix}
\vecnull\\
\vecalf
\end{pmatrix})
\overline{\Theta_g(u+\i\frac1\lambda,0;\begin{pmatrix}
\vecnull\\
\vecalf
\end{pmatrix})}
$$
with
$$
f(\vecw)=\psi_1(\|\vecw\|^2),\qquad g(\vecw)=\psi_2(\|\vecw\|^2).
$$

\section{Equidistribution\label{secequi}}

\subsection{Theorem\label{equithm}}{\em
Let $\Gamma$ be a subgroup of $\SL(2,\ZZ)\ltimes\ZZ^{2k}$ of finite index,
and assume the components of the vector $(\vecy,1)\in\RR^{k+1}$ are linearly
independent over $\QQ$.
Let $h$ be a continuous function $\RR\rightarrow\RR_+$
with compact support.
Then, for any bounded continuous function $F$ on 
$\Gamma\backslash G^k$ and any $\sigma\geq 0$, we have 
$$
\lim_{v\rightarrow 0} 
v^\sigma \int_\RR F(u+\i v,0; \begin{pmatrix}
\vecnull\\
\vecy
\end{pmatrix}) 
\; h(v^\sigma u)\, du
= \frac{1}{\mu(\Gamma\backslash G^k)}\int_{\Gamma\backslash G^k} F \, d\mu
\;\int h(w)\, dw
$$
where $\mu$ is the Haar measure of $G^{k}$.
}

\begin{proof}
For $\sigma=0$ the above statement is proved in \cite{Marklof00},
Theorem 5.7; see also Shah's more general Theorem 1.4 in \cite{Shah96}. 
The case $\sigma>0$ is easier and in fact follows from the result
for $\sigma=0$, since the translate of the unipotent
orbit is expanding at a faster rate:

As in \cite{Marklof00}, Section 5, we define the unipotent flow
$\Psi^t: \Gamma\backslash G^k\rightarrow\Gamma\backslash G^k$ 
by right translation with
$$
\Psi_0^t=
(\begin{pmatrix}
1 & t \\
0 & 1
\end{pmatrix} ; \vecnull ),
$$
and furthermore the flow 
$\Phi^t: \Gamma\backslash G^k\rightarrow\Gamma\backslash G^k$ 
by right translation with
$$
\Phi_0^t=
(\begin{pmatrix}
\e^{-t/2} & 0 \\
0 & \e^{t/2}
\end{pmatrix} ; \vecnull ).
$$
By Theorem 5.7 in \cite{Marklof00}, the orbit segment
$$
\Gamma \{ (u+\i \e^{-t},0; \begin{pmatrix}
\vecnull\\
\vecy
\end{pmatrix}) 
: \, u\in[-1,1]
\}
$$
is dense in $\Gamma\backslash G^k$
in the limit $t\rightarrow\infty$. Hence we find a sequence
$\{ u_t \}_{t\in\RR_+}$ with $u_t\in[-1,1]$ such that
$$
\Gamma g_t :=\Gamma(u_t+\i \e^{-t},0; \begin{pmatrix}
\vecnull\\
\vecy
\end{pmatrix}) 
=
\Gamma (1;\begin{pmatrix} \vecnull \\ \vecy \end{pmatrix})
\Psi^{u_t}\Phi^t
$$
converges in the limit $t\rightarrow\infty$ 
to a generic point in $\Gamma\backslash G^k$. 
Theorem 2 in \cite{Dani93} implies then that
for any constant $B\neq 0$
\begin{multline*}
\frac{1}{B\e^{\sigma t}} \int_0^{B\e^{ \sigma t}} F(u_t+u+\i\e^{-t},0; 
\begin{pmatrix}
\vecnull\\
\vecy
\end{pmatrix}) \, du 
=
\frac{1}{B\e^{(1+\sigma)t}} \int_0^{B\e^{(1+\sigma)t}} F(g_t \Psi^u) \, du  \\
\rightarrow
\frac{1}{\mu(\Gamma\backslash G^k)}\int_{\Gamma\backslash G^k} F \, d\mu 
\end{multline*}
as $t\rightarrow\infty$. 
Because $F$ is bounded and $u_t$ is contained in a compact interval, 
note that
\begin{multline*}
\frac{1}{B \e^{\sigma t}} \int_0^{B\e^{\sigma t}} F(u_t+u+\i\e^{-t},0; 
\begin{pmatrix}
\vecnull\\
\vecy
\end{pmatrix}) \, du 
=
\frac{1}{B \e^{\sigma t}} \int_{u_t}^{B \e^{\sigma t}+u_t} F(u+\i\e^{-t},0; 
\begin{pmatrix}
\vecnull\\
\vecy
\end{pmatrix}) \, du \\
=
\frac{1}{B \e^{\sigma t}} \int_0^{B\e^{\sigma t}} F(u+\i\e^{-t},0; 
\begin{pmatrix}
\vecnull\\
\vecy
\end{pmatrix}) \, du +O(\e^{-\sigma t}).
\end{multline*}
Therefore, for any constants $-\infty<A<B<\infty$,
$$
\lim_{t\rightarrow\infty}
\frac{1}{\e^{\sigma t}} 
\int_{A \e^{\sigma t}}^{B \e^{\sigma t}}  
F(u+\i\e^{-t},0; \begin{pmatrix}
\vecnull\\
\vecy
\end{pmatrix}) \, du =
\frac{(B-A)}{\mu(\Gamma\backslash G^k)}\int_{\Gamma\backslash G^k} F \, d\mu .
$$
The theorem now follows from a standard approximation argument
(approximate $h$ from above and below by step functions).
\end{proof}

\section{Diophantine conditions\label{secdiophantine}}

\subsection{}

In order to extend the equidistribution results
to unbounded test functions
such as $\Theta_f\overline\Theta_g$, let us
study the following model functions, whose asymptotics in the
cusp is similar to that of $\Theta_f\overline\Theta_g$.
Let $G=G^k$ and $\Gamma=\SL(2,\ZZ)\ltimes \ZZ^{2k}$. Define
furthermore the subgroup
$$
\Gamma_\infty = \{ \begin{pmatrix} 1 & m \\ 0 & 1 \end{pmatrix}:
m\in\ZZ\} \subset \SL(2,\ZZ),
$$
and put 
$$
v_\gamma := \Im (\gamma\tau) = \frac{v}{|c\tau+d|^2}, \quad
\text{ for } \gamma=\begin{pmatrix} a & b \\ c & d \end{pmatrix},
$$
and
$\vecy_\gamma := c \vecx + d \vecy$.
Let $\chi_R$ be the characteristic function of the interval $[R,\infty)$,
$$
\chi_R(t)=
\begin{cases} 
1 & (t\geq R) \\
0 & (t < R) .
\end{cases}
$$
For any $f\in\C(\RR^k)$ of rapid decay
(i.e., $f(\vecw)$ decays rapidly for $\|w\|\rightarrow\infty$) the function
$$
F_R(\tau;\vecxi)= 
\sum_{\gamma\in\Gamma_\infty\backslash\SL(2,\ZZ)}
\sum_{\vecm\in\ZZ^k} 
f\big( (\vecy_\gamma+\vecm) v_\gamma^{1/2} \big)\,
v_\gamma^\beta \, \chi_R(v_\gamma), \quad R>1,
$$
is invariant under the action of $\Gamma$.
If $\tau$ lies in the fundamental domain of $\SL(2,\ZZ)$,
given by 
${\mathcal F}_{\SL(2,\ZZ)}
=\{ \tau\in\H: u\in[-\tfrac12,\tfrac12), |\tau|> 1 \}$,
then $F_R(\tau;\vecxi)$ has the representation
$$
F_R(\tau;\vecxi)= 
\sum_{\vecm\in\ZZ^k} 
\big\{ f\big( (\vecy+\vecm) v^{1/2} \big) 
+f\big((-\vecy+\vecm) v^{1/2} \big) \big\}
v^\beta \chi_R(v) .
$$
The remaining sum over $\vecm$ is rapidly converging since $f$
is of rapid decay.

\subsection{\label{lone}}

The $\L^1$ norm of $F_R$ over $\Gamma\backslash G$ is, for $f\geq 0$,
$$
\mu(F_R)=\int_{\Gamma\backslash G} F_R(\tau;\vecxi)\; d\mu(\tau,\phi;\vecxi)
$$
with Haar measure 
$$
d\mu(\tau,\phi;\vecxi)=\frac{du\,dv\,d\phi\,dx\,dy}{v^2}.
$$
We therefore have
$$
\mu(F_R)=2\pi \int_{\RR^k}  f(\vecw) dw  \int_R^\infty v^{\beta-k/2-2} dv
= 2\pi \frac{R^{-(k/2+1-\beta)}}{k/2+1-\beta}  \int_{\RR^k}  f(\vecw) dw
$$
for $\beta<k/2+1$, and $\mu(F_R)=\infty$ otherwise.
In the following we will be especially interested in $\beta=k/2$, for which
$$
\mu(F_R)= 2\pi R^{-1} \int_{\RR^k}  f(\vecw) dw .
$$

\subsection{\label{sun}}

As in Section 6.4. \cite{Marklof00} 
we may write the sum in $F_R(\tau;\vecxi)$ explicitly as
\begin{multline*}
F_R(\tau;\vecxi)=
\sum_{\vecm\in\ZZ^k} \big\{ f\big( (\vecy+\vecm) v^{1/2} \big) 
+ f\big( (-\vecy+\vecm) v^{1/2} \big) \big\}
v^\beta \chi_R(v) \\
+
\sum_{\vecm\in\ZZ^k} 
\big\{ f\big((\vecx + \vecm) \frac{v^{1/2}}{|\tau|} \big) +
f\big((-\vecx + \vecm) \frac{v^{1/2}}{|\tau|} \big) \big\}
\frac{v^\beta}{|\tau|^{2\beta}} \chi_R(\frac{v}{|\tau|^2}) \\
+
\sum_{\substack{(c,d)\in\ZZ^2\\ \gcd(c,d)=1 \\ c,d\neq 0}}
\sum_{\vecm\in\ZZ^k} 
f\big((c \vecx + d \vecy + \vecm) \frac{v^{1/2}}{|c\tau+d|} \big)
\frac{v^\beta}{|c\tau+d|^{2\beta}} \chi_R(\frac{v}{|c\tau+d|^2}).
\end{multline*}

In what follows we will restrict our attention to the case 
$\beta=k/2$ and $\vecxi=(\begin{smallmatrix}
\vecnull\\
\vecy
\end{smallmatrix})$.

\subsection{Proposition\label{rain}}{\em
Let $\vecy$ be diophantine of type $\kappa$.
Then, for any $\epsilon,\epsilon'$ with $0<\epsilon<1$ and
$0<\epsilon'<\frac{1}{\kappa-1}$, 
$$
\limsup_{v\rightarrow 0}
v^{k/2-1}
\int_{|u|> v^{1-\epsilon}} F_R(u+\i v;\begin{pmatrix}
\vecnull\\
\vecy
\end{pmatrix}) \; h(v^{k/2-1}u)\, du 
\ll_{\epsilon,\epsilon'}  R^{-(\frac{1}{\kappa-1}-k+2)/2} + R^{-\epsilon'/2}.
$$
}

Note that the above expression
vanishes, for $R\rightarrow\infty$, when $\kappa<\frac{k-1}{k-2}$.
The second term is obviously only relevant in dimension $k=2$,
since for $k>2$ we may chose $\epsilon'$ in such a way that
$\frac{1}{\kappa-1}<\epsilon'+k-2$.

The key ingredient in the proof is the following lemma.

\subsection{Lemma\label{moon}} {\em
Let $\vecalf$ be diophantine of type $\kappa$, and $f\in\C(\RR^k)$
of rapid decay.
Then, for any fixed $A>1$ and $\epsilon>0$ with $\epsilon<\frac{1}{\kappa-1}$,
$$
\sum_{d=1}^D \sum_{\vecm\in\ZZ^k} 
f\big( T (d\vecalf+\vecm) \big) 
\ll
\begin{cases}
T^{-A} & (D \leq T^\epsilon) \\[10pt]
1 & (T^\epsilon \leq D \leq T^{\frac{1}{\kappa-1}}) \\[10pt]
D\, T^{-\frac{1}{\kappa-1}}   & (D \geq T^{\frac{1}{\kappa-1}}) ,
\end{cases}
$$
uniformly for all $D,T>1$.} 

\subsection{\it Proof.}
Let us divide the sum over $d$ into blocks of the form
$$
\sum_{0\leq d \leq T^{\frac{1}{\kappa-1}}}  
\sum_{\vecm\in\ZZ^k} 
f\big( T ((b+d)\vecalf+\vecm) \big) .
$$
The number of such blocks is $\ll D T^{-\frac{1}{\kappa-1}} + 1$.
Since $\vecalf$ is of type $\kappa$ there is a constant $C$ such that,
for all $0< |q| \leq T^{\frac{1}{\kappa-1}}$
we have 
$$
\frac{C}{|q| T} \leq \frac{C}{|q|^{\kappa}} 
\leq \max_j |\alpha_j -\frac{m_j}{q}|,
$$
and thus
$$
\max_j |q\alpha_j - m_j| \geq \frac{C}{T} .
$$
For $b$ fixed, the minimal distance between the points
$(b+d)\vecalf+\vecm$ ($0\leq d \leq T^{\frac{1}{\kappa-1}}$, $\vecm\in\ZZ^k$)
is bounded from below by
$$
\min_{0<|q| \leq T^{\frac{1}{\kappa-1}},\vecm\in\ZZ^k}
\| q\vecalf+\vecm \| \geq 
\min_{0<|q| \leq T^{\frac{1}{\kappa-1}},\vecm\in\ZZ^k} \max_j |q\alpha_j - m_j|
\geq \frac{C}{T} .
$$
Hence any rectangular box with sides $\ll \frac1T$ contains at most
a bounded number of points. Because $f$ is rapidly decreasing,
we therefore find 
$$
\sum_{0\leq d \leq T^{\frac{1}{\kappa-1}}} 
\sum_{\vecm\in\ZZ^k} 
f\big( T ((b+d)\vecalf+\vecm) \big) \ll 1 ,
$$
independently of $b$. This explains the second and third bound.
The first bound is obtained from
$$
\| d \vecalf+\vecm \| \geq \max_j |d\alpha_j - m_j|
\geq \frac{C}{d^{\kappa-1}} \geq \frac{C}{D^{\kappa-1}}
$$
which holds for all $d=1,\ldots, D$. Since $f$ is rapidly decreasing, we
have
$$
\sum_{d=1}^D \sum_{\vecm\in\ZZ^k} 
f\big( T (d\vecalf+\vecm) \big) 
\ll D (\frac{D^{\kappa-1}}{T})^B
$$
for any $B>1$.
\qed

\subsection{Proof of Proposition \ref{rain}}

Let us assume without loss of generality that $f$ is positive and even, 
i.e., $f\geq 0$, $f(-\vecw)=f(\vecw)$.

It follows from the expansion in \ref{sun} that, for $v<1$, 
the first term involving $\chi_R(v)$ vanishes and hence we are left with
\begin{multline*}
F_R(\tau;\begin{pmatrix}
\vecnull\\
\vecy
\end{pmatrix})= 
2 \sum_{\vecm\in\ZZ^k} 
f\big(\vecm \frac{v^{1/2}}{|\tau|} \big) 
\frac{v^{k/2}}{|\tau|^k} \chi_R(\frac{v}{|\tau|^2}) \\
+2\sum_{\substack{(c,d)\in\ZZ^2\\ \gcd(c,d)=1 \\ c>0,d\neq 0}}
\sum_{\vecm\in\ZZ^k} 
f\big((d \vecy + \vecm) \frac{v^{1/2}}{|c\tau+d|} \big)
\frac{v^{k/2}}{|c\tau+d|^k} \chi_R(\frac{v}{|c\tau+d|^2}).
\end{multline*}

\subsubsection{}
With regard to the first term in the above expansion, a change of variable
$u=vt$ yields
\begin{multline*}
v^{k/2-1}
\int_{|u|>v^{1-\epsilon}}  2 \sum_{\vecm\in\ZZ^k} 
f\big(\vecm \frac{v^{1/2}}{|\tau|} \big) 
\frac{v^{k/2}}{|\tau|^k} \chi_R(\frac{v}{|\tau|^2}) h(v^{k/2-1}u)
\, du \\
= 2 \sum_{\vecm\in\ZZ^k} \int_{|t|>v^{-\epsilon}}
f\big(\frac{\vecm}{v^{1/2} (t^2+1)^{1/2}} \big) 
\frac{1}{(t^2+1)^{k/2}} 
\chi_R(\frac{1}{v(t^2+1)}) h(v^{k/2} t)\, dt \\
\sim 2 f(0)\, h(0) \int_{|t|>v^{-\epsilon}} \frac{dt}{(t^2+1)^{k/2}}
\rightarrow 0, 
\end{multline*}
as $v\rightarrow 0$.

\subsubsection{}
An upper bound for the remaining terms is obtained by dropping
the condition $|u|>v^{1-\epsilon}$ in the integral. We then need
to estimate
$$
S(v)=
\sum_{\substack{(c,d)\in\ZZ^2\\ \gcd(c,d)=1 \\ c>0,d\neq 0}}
\sum_{\vecm\in\ZZ^k} 
J(v,c,d,\vecm)
$$
with
$$
J(v,c,d,\vecm) = v^{k/2-1}
\int_\RR f\big((d \vecy + \vecm) \frac{v^{1/2}}{|c\tau+d|} \big)
\frac{v^{k/2}}{|c\tau+d|^k}  \chi_R(\frac{v}{|c\tau+d|^2}) \; 
h(v^{k/2-1}u)\, du.
$$
Substituting $u$ by $t=v^{-1}(u+\frac{d}{c})$ gives
$$
\frac{1}{c^k}
\int_\RR f\big((d \vecy + \vecm) \frac{1}{\sqrt{c^2 v(t^2+1)}} \big)
\frac{1}{(t^2+1)^{k/2}} \chi_R(\frac{1}{c^2v(t^2+1)}) \; h\big(v^{k/2-1}
(vt-\frac{d}{c})\big)\, dt.
$$
The range of integration is bounded by
$$
R < \frac{1}{c^2 v(t^2+1)}, \qquad\text{i.e.}\quad
|t| \ll \frac{1}{c \sqrt{vR}}.
$$
Therefore $|vt|\ll v^{1/2}c^{-1} R^{-1/2}$ is 
uniformly close to zero, and hence, because of the compact support of $h$,
we find $|d|\ll cv^{-(k/2-1)}$. So
$$
S(v) \ll
\sum_{c=1}^\infty
\sum_{0<|d|\ll c v^{-(k/2-1)}}
\sum_{\vecm\in\ZZ^k} 
K(v,c,d,\vecm) ,
$$
with 
$$
K(v,c,d,\vecm)=
\frac{1}{c^k}
\int_\RR f\big((d \vecy + \vecm) \frac{1}{\sqrt{c^2 v(t^2+1)}} \big)
\frac{1}{(t^2+1)^{k/2}} \chi_R(\frac{1}{c^2v(t^2+1)}) \; dt.
$$

\subsubsection{}
To apply Lemma \ref{moon} with $D=cv^{-(k/2-1)}$,
$T=(c^2 v(t^2+1))^{-1/2}>\sqrt{R}>1$, split the range of integration
into the ranges
\begin{equation*}
\begin{split}
(1): & \quad c v^{-(k/2-1)}\leq (c^2 v(t^2+1))^{-\epsilon/2} \\
(2): & \quad (c^2 v(t^2+1))^{-\epsilon/2} \leq c v^{-(k/2-1)}
\leq (c^2 v(t^2+1))^{-\delta/2} \\
(3): & \quad c v^{-(k/2-1)} \geq (c^2 v(t^2+1))^{-\delta/2}
\end{split}
\end{equation*}
with $\delta=\frac{1}{\kappa-1}$.
Denote the corresponding integrals by
$K_1(v,c,d,\vecm)$, $K_2(v,c,d,\vecm)$ and $K_3(v,c,d,\vecm)$, respectively.

\subsubsection{}
Because $R^{-1/2}\geq T^{-1}$, we find in the first range
$D\leq T^\epsilon$ that 
\begin{equation*}
\begin{split}
\sum_{c>0} \sum_{d\ll c v^{-(k/2-1)}} 
\sum_{\vecm\in\ZZ^k} 
K_1(v,c,d,\vecm) 
& \ll  R^{-A/2} \sum_{c>0}\frac{1}{c^k}
\int_{(1)}   \frac{1}{(t^2+1)^{k/2}} \chi_R(\frac{1}{c^2v(t^2+1)}) \;  dt \\
& \ll  R^{-A/2} \sum_{c>0}\frac{1}{c^k}
\int_\RR   \frac{1}{(t^2+1)^{k/2}} \;  dt \\
& \ll R^{-A/2} .
\end{split}
\end{equation*}

\subsubsection{}
For an upper bound
the second range $T^\epsilon\leq D \leq T^\delta$ may be extended to
$T^{\epsilon}\leq D$, i.e., 
$$
c^{1+\epsilon} (t^2+1)^{\epsilon/2} \geq v^{k/2-1-\epsilon/2} .
$$ 
We have therefore
\begin{equation*}
\begin{split}
\sum_{c>0} \sum_{d\ll c v^{-(k/2-1)}} 
\sum_{\vecm\in\ZZ^k} 
K_2(v,c,d,\vecm) 
& \ll  \sum_{c>0} \frac{1}{c^k}
\int_{(2)}   \frac{dt}{(t^2+1)^{k/2}} \\
& \ll  \sum_{c>0} \frac{1}{c^k}
\big\{ c^{1+\epsilon} v^{-(k/2-1-\epsilon/2)} 
\big\}^{(k/2-1)2/\epsilon} 
\int_{\RR}   \frac{dt}{t^2+1} \\
& \ll v^A
\sum_{c>0} c^{-B} 
\end{split}
\end{equation*}
with
$$
A=-\big(\frac{k}{2}-1-\frac{\epsilon}{2}\big)\big(\frac{k}{2}-1\big)
\frac{2}{\epsilon}
$$ 
and 
$$
B=-\big(\frac{k}{2}- 1 - \epsilon\big)\frac{2}{\epsilon}
= 1 -\big(\frac{k}{2}- 1 - \frac{\epsilon}{2}\big)\frac{2}{\epsilon}.
$$
If we chose $\epsilon$ in a way that $k-2<\epsilon<\delta=\frac{1}{\kappa-1}$,
we find that for $k>2$ we have $A>0$ and $B>1$.
Hence
$$
\sum_{c>0} \sum_{d\ll c v^{-(k/2-1)}} 
\sum_{\vecm\in\ZZ^k}  
K_2(v,c,d,\vecm) \rightarrow 0 
$$
for small $v$.
In the case $k=2$ we exploit 
the inclusion $R^{\epsilon/2}<T^{\epsilon}\leq D\ll c$, which yields
$$
\sum_{c>0} \sum_{d\ll c} \sum_{\vecm\in\ZZ^2} K_2(v,c,d,\vecm) 
\ll  \sum_{c>R^{\epsilon/2}} c^{-2}
\int   \frac{dt}{t^2+1}  
\ll R^{-\epsilon/2} ,
$$
compare \cite{Marklof00}.

\subsubsection{}
In the third range, we have for $v$ sufficiently small
\begin{equation*}
\begin{split}
\sum_{c>0}& \sum_{d\ll cv^{-(k/2-1)}} \sum_{\vecm\in\ZZ^k} K_3(v,c,d,\vecm) \\
& \ll  \sum_{c>0} \frac{1}{c^k} cv^{-(k/2-1)}
\int_{(3)} c^\delta v^{\delta/2} (t^2+1)^{(\delta-k)/2} 
\chi_R(\frac{1}{c^2v(t^2+1)}) \; dt \\ 
& = v^{(\delta-k)/2+1} \sum_{c>0} c^{\delta-k+1} 
\int_{(3)}  
(t^2+1)^{(\delta-k)/2} 
\chi_R(\frac{1}{c^2v(t^2+1)}) \;  dt  \\
& \ll v^{(\delta-k)/2+1} \int_\RR \big\{ \sum_{c=1}^\infty c^{\delta-k+1}  
\chi_R(\frac{1}{c^2v(t^2+1)}) \big\} (t^2+1)^{(\delta-k)/2}  \, dt \\
& < v^{(\delta-k)/2+1} \int_\RR \big\{ \int_0^\infty x^{\delta-k+1} 
\chi_R(\frac{1}{x^2v(t^2+1)}) dx \big\} (t^2+1)^{(\delta-k)/2}  \, dt\\  
& = \int_\RR \big\{ \int_0^\infty  x^{\delta-k+1}
\chi_R(\frac{1}{x^2}) dx \big\} (t^2+1)^{-1} \, dt  \\
& = \int_\RR \big\{ \frac{x^{\delta-k+2}}{\delta-k+2} \big\}_0^{R^{-1/2}} 
(t^2+1)^{-1} \, dt  \\
& = \pi \frac{R^{-(\delta-k+2)/2}}{\delta-k+2}   .
\end{split}
\end{equation*}
The proof of Proposition \ref{rain} is complete.
\qed
\\

\subsection{\label{FR}}
Let us define the characteristic function on $\Gamma\backslash G^k$
$$
X_R(\tau)= 
\sum_{\gamma\in\{\Gamma_\infty\cup(-1)\Gamma_\infty\}\backslash\SL(2,\ZZ)}
\chi_R(v_\gamma),
$$
where $\chi_R$ is the characteristic function of $[R,\infty)$.
Proposition \ref{rain} allows us now to extend the equidistribution
theorem \ref{equithm} to unbounded functions which are dominated
by $F_R$, i.e. that is, for some fixed constant $L>1$ we have
$$
|F(\tau,\phi;\vecxi)|X_R(\tau) \leq L+F_R(\tau;\vecxi)
$$
for all sufficiently
large $R>1$, 
uniformly for all $(\tau,\phi;\vecxi)\in G^k$.

\subsection{Theorem\label{propo}}{\em
Let $\Gamma$ be a subgroup of $\SL(2,\ZZ)\ltimes\ZZ^{2k}$ of finite index.
Let $h$ be a continuous function $\RR\rightarrow\RR_+$
with compact support. 
Suppose the continuous function $F\geq 0$ is dominated by $F_R$.
Fix some $\vecy\in\TT^{k}$
such that the components of the vector 
$(\vecy,1)\in\RR^{k+1}$ are linearly
independent over $\QQ$.
Then, for any $\epsilon$ with $0<\epsilon<1$,
$$
\liminf_{v\rightarrow 0}
v^{k/2-1}\int_{|u|> v^{1-\epsilon}} F(u+\i v,0;\begin{pmatrix}
\vecnull\\
\vecy
\end{pmatrix}) 
\; h(v^{k/2-1}u)\, du 
\geq \frac{1}{\mu(\Gamma\backslash G^k)}\int_{\Gamma\backslash G^k} F \, d\mu 
\; \int h .
$$
Assume furthermore
that $\vecy$ is diophantine of type $\kappa<\frac{k-1}{k-2}$.
Then, for any $\epsilon>0$,
$$
\limsup_{v\rightarrow 0}
v^{k/2-1}\int_{|u|> v^{1-\epsilon}} F(u+\i v,0;\begin{pmatrix}
\vecnull\\
\vecy
\end{pmatrix}) 
\; h(v^{k/2-1}u)\, du 
\leq \frac{1}{\mu(\Gamma\backslash G^k)}\int_{\Gamma\backslash G^k} F \, d\mu 
\; \int h.
$$
}

\begin{proof}
The theorem follows from Theorem \ref{equithm} and Proposition \ref{rain}
in the identical manner as Theorem 7.3 in \cite{Marklof00}.
\end{proof}
 
\subsection{}
The subgroup $\Gamma=\Gamma^k$ is a subgroup
of finite index in $\SL(2,\ZZ)\ltimes (\tfrac12\ZZ)^{2k}$
rather than $\SL(2,\ZZ)\ltimes\ZZ^{2k}$ (Lemma \ref{two}).
We therefore need to rephrase Theorem \ref{propo} slightly.
Define the dominating function $\hat F_R$ on 
$\Gamma\backslash G^k$ by 
$\hat F_R(\tau;\vecxi)=F_R(\tau; 2 \vecxi)$,
with $F_R$ as in \ref{FR}.

\subsection{Corollary}\label{cordio1}
{\em 
Let $\Gamma$ be a subgroup of $\SL(2,\ZZ)\ltimes (\tfrac12\ZZ)^{2k}$
of finite index, 
$h$, $\vecy$ as in Theorem \ref{propo},
and $F: \Gamma\backslash G^k\rightarrow\CC$
a continuous function which is dominated by $\hat F_R$.
If $\vecy$ is diophantine of type $\kappa<\frac{k-1}{k-2}$.
Then, for any $\epsilon$ with $0<\epsilon<1$,
$$
\lim_{v\rightarrow 0}
v^{k/2-1}\int_{|u|> v^{1-\epsilon}} F(u+\i v,0;\begin{pmatrix}
\vecnull\\
\vecy
\end{pmatrix}) 
\; h(v^{k/2-1}u)\, du 
=\frac{1}{\mu(\Gamma\backslash G^k)}\int_{\Gamma\backslash G^k} F \, d\mu 
\; \int h.
$$
}

\begin{proof}
The proof is analogous to that of Corollary 7.6 in \cite{Marklof00}.
\end{proof}

\section{The Main Theorem\label{secmain}}

\subsection{Main Theorem\label{mthm}}{\em
Suppose $f(\vecw)=\psi_1(\| \vecw \|^2)$ and $g(\vecw)
=\psi_2(\| \vecw \|^2)$ with $\psi_1,\psi_2\in\scrs(\RR_+)$ real-valued.
Let $h$ be a continuous function $\RR\rightarrow\CC$
with compact support. 
Assume that the components of $(\vecy,1)\in\RR^{k+1}$ 
are linearly independent over $\QQ$ and that
$\vecy$ is diophantine of type $\kappa<\frac{k-1}{k-2}$.
Then
\begin{multline*}
\lim_{v\rightarrow 0}
v^{k/2-1}\int_\RR \Theta_f(u+\i v,0;\begin{pmatrix}
\vecnull\\
\vecy
\end{pmatrix})
\overline{\Theta_g(u+\i v,0;\begin{pmatrix}
\vecnull\\
\vecy
\end{pmatrix})} \; h(v^{k/2-1} u)\, du \\
= \frac{k^2}{2}B_k^2 h(0) \int_0^\infty \psi_1(r)\psi_2(r) r^{k-2}\, dr \\
+ \frac{k}{2}B_k \int_\RR h(u)\, du 
\int_0^\infty \psi_1(r)\psi_2(r) r^{k/2-1}\, dr , 
\end{multline*}
where $B_k$ is the volume of the $k$-dimensional unit ball.} \\
 
We will need the following two lemmas.

\subsection{Lemma\label{mlemm1}}{\em We have
$$
\frac{1}{\mu(\Gamma^k\backslash G^k)}
\int_{\Gamma^k\backslash G^k} 
\Theta_f(\tau,\phi;\vecxi)\overline{\Theta_g(\tau,\phi;\vecxi)}\, d\mu
= \int_{\RR^k} f(\vecw)\overline{g(\vecw)} dw .
$$
}

Note that if $f(\vecw)=\psi_1(\| \vecw \|^2)$ and $g(\vecw)
=\psi_2(\| \vecw \|^2)$, we have
$$
\int f(\vecw)\overline{g(\vecw)} \, dw
= \frac{k}{2}B_k \int_0^\infty \psi_1(r)\psi_2(r) r^{k/2-1}\, dr .
$$ 

\begin{proof}
A short calculation shows that
$$
\int_{\TT^{2k}}
\Theta_f(\tau,\phi;\vecxi)\overline{\Theta_g(\tau,\phi;\vecxi)} \, d\xi  
= \int f_\phi(\vecw)\overline{g_\phi(\vecw)} \, dw .
$$
Since $f_\phi=U^\phi f$ with $U^\phi$ unitary, we have
$$
\int f_\phi(\vecw)\overline{g_\phi(\vecw)} \, dw =
\int f(\vecw)\overline{g(\vecw)} \, dw .
$$
\end{proof}

\subsection{Lemma\label{mlemm2}}{\em 
Suppose $f(\vecw)=\psi_1(\| \vecw \|^2)$ and $g(\vecw)
=\psi_2(\| \vecw \|^2)$.
For any $\tfrac12<\gamma<1$, we have
\begin{multline*}
\lim_{v\rightarrow 0} v^{k/2-1} \int_{|u|< v^{\gamma}}
\Theta_f(u+\i v, 0 ;\begin{pmatrix}
\vecnull\\
\vecy
\end{pmatrix})
\overline{\Theta_g(u+\i v, 0 ;\begin{pmatrix}
\vecnull\\
\vecy
\end{pmatrix})} h(v^{k/2-1} u)\,du  \\
= \frac{k^2}{2} B_k^2 h(0) \int_0^\infty \psi_1(r)\psi_2(r)\, r^{k-2}dr .
\end{multline*}
}

\begin{proof}
From Proposition \ref{asymp} we know that
$$
\Theta_f(-\frac{1}{\tau},\arg\tau;\begin{pmatrix}
-\vecy\\
\vecnull
\end{pmatrix})
\overline{\Theta_g(-\frac{1}{\tau},\arg\tau;\begin{pmatrix}
-\vecy\\
\vecnull
\end{pmatrix})} 
= \frac{v^{k/2}}{|\tau|^k} f_{\arg\tau}(\vecnull) 
\overline{g_{\arg\tau}(\vecnull)}
+O_R((\frac{v}{|\tau|^2})^{-R})
$$
holds uniformly for $|u|<v^{1/2}<1$. The remainder vanishes 
for $|u|<v^{\gamma}<1$.
Now 
\begin{multline*}
f_{\arg\tau}(\vecnull)\overline{g_{\arg\tau}(\vecnull)} 
= \frac{|\tau|^k}{v^k}
\{\int e(\tfrac12 \|\vecw\|^2 \frac{u}{v})\, f(\vecw)\, dw \}
\overline{\{\int e(\tfrac12 \|\vecw\|^2 \frac{u}{v})\, g(\vecw)\, dw\}}\\
= \frac{|\tau|^k}{v^k}\, \frac{k^2}{4} B_k^2 
\int_0^\infty e(\frac{(r_1-r_2)u}{2v}) \psi_1(r_1)\psi_2(r_2)\, 
r_1^{k/2-1}dr_1\, r_2^{k/2-1} dr_2
\end{multline*}
(substitute $\vecw$ by polar coordinates) and so, as $v\rightarrow\infty$,
\begin{equation*}
\begin{split}
\int_{|u|< v^{\gamma}} &
v^{k/2-1} \Theta_f(u+\i v,0;\begin{pmatrix}
\vecnull\\
\vecy
\end{pmatrix})
\overline{\Theta_g(u+\i v,0;\begin{pmatrix}
\vecnull\\
\vecy
\end{pmatrix})} h(v^{k/2-1} u)\,du \\
& \sim v^{-1} \frac{k^2}{4} B_k^2 \int_{|u|< v^{\gamma}}
\int_0^\infty e(\frac{(r_1-r_2)u}{2v}) \psi_1(r_1)\psi_2(r_2)\,
\\ & \qquad\qquad\qquad\qquad  
r_1^{k/2-1} dr_1\,  r_2^{k/2-1} dr_2
 h(v^{k/2-1}u)\,du \\
& \sim \frac{k^2}{2} B_k^2  h(0) \int_{2|u|< v^{\gamma-1}}
\int_0^\infty e((r_1-r_2)u) \\ & 
\qquad\qquad\qquad\qquad \psi_1(r_1)\psi_2(r_2)\, 
r_1^{k/2-1} dr_1\,  r_2^{k/2-1} dr_2 \, du \\
& \sim  \frac{k^2}{2} B_k^2  h(0) \int_\RR
\int_0^\infty e((r_1-r_2)u) \psi_1(r_1)\psi_2(r_2)\, 
\\ & \qquad\qquad\qquad\qquad 
r_1^{k/2-1} dr_1\,  r_2^{k/2-1} dr_2 \, du  \\
& = \frac{k^2}{2} B_k^2 h(0) \int_0^\infty \psi_1(r)\psi_2(r)\, r^{k-2}dr 
\end{split}
\end{equation*}
by Parseval's equality.
\end{proof}

\subsection{Proof of the Main Theorem}

We may assume without 
loss of generality that in Theorem \ref{mthm} $h$ is positive.
Split the integration on the left-hand-side of
\ref{mthm} into 
$$
\int_\RR = \int_{|u|< v^{1-\epsilon}} + \int_{|u|> v^{1-\epsilon}} ,
$$
for some small $\epsilon>0$.
The first integral gives, by virtue of Lemma \ref{mlemm2},
the contribution
$$
\frac{k^2}{2} B_k^2 h(0) \int_0^\infty \psi_1(r)\psi_2(r)\, r^{k-2}dr
$$
Corollary \ref{cordio1}, together with Lemma \ref{mlemm1},
yields the second term on the right-hand-side of \ref{mthm}.
Compare Section 8.4 in \cite{Marklof00} for more details.
\qed

\subsection{Proof of Theorem \ref{lthm}}

We have by construction
$$
R_2(\psi_1,\psi_2,h,\lambda) =
\frac{1}{B_k}
v^{k/2-1}\int_\RR 
\Theta_f(u+\i \frac{1}{\lambda},0;\begin{pmatrix}
\vecnull\\
\vecalf
\end{pmatrix})
\overline{\Theta_g(u+\i \frac{1}{\lambda},0;\begin{pmatrix}
\vecnull\\
\vecalf
\end{pmatrix})} \; 
h(v^{k/2-1}u)\, du
$$
with $v=\lambda^{-1}$. 
Recall that $2h(0)=\int \hat h(s)ds$ by Fourier inversion
and thus we have finally $\int h(u) du = \hat h(0)$.
\qed

\section{Counter examples\label{seccounter}}

We assume throughout this section that $k>2$. The case 
$k=2$ is studied in  \cite{Marklof00}, Section 9.

\subsection{\label{bubu}}
Suppose $\alpha_{k-1},\alpha_k$ 
are both rational and
$(\alpha_1,\ldots,\alpha_{k-2})$ is a badly approximable $(k-2)$-tuple.
In this case we find a constant $C$ such that
$$
\max_{1\leq j\leq k} |\alpha_j - \frac{m_j}{q}| \geq
\max_{1\leq j\leq k-2} |\alpha_j - \frac{m_j}{q}| 
> \frac{C}{q^{1+\frac{1}{k-2}}}
$$
for all $m_1,\ldots,m_j,q\in\ZZ$, $q>0$,
and so $\vecalf$ is of type $\kappa=\frac{k-1}{k-2}$.
On the other hand we have
\begin{equation*}
\begin{split}
\# \{ (\vecm,\vecn)& \in\ZZ^k\times\ZZ^k: \; \vecm\neq\vecn,\; \\ 
&\|\vecm-\vecalf\|^k \leq X, \|\vecn-\vecalf\|^k \leq X, \; 
 \|\vecm-\vecalf\|^k=\|\vecn-\vecalf\|^k \}
\\
\geq
\# \{ (\vecm,\vecn)&\in\ZZ^k\times\ZZ^k:  \; \vecm\neq\vecn,\; 
(m_1,\ldots,m_{k-2})=(n_1,\ldots,n_{k-2}), \\
& \|\vecm-\vecalf\|^k \leq X , \|\vecn-\vecalf\|^k \leq X, \; 
\|\vecm-\vecalf\|^2=\|\vecn-\vecalf\|^2 \} .
\end{split}
\end{equation*}
This is easily seen to be bounded from below by
\begin{multline*}
\gg X^{(k-2)/k} \#\{ (m_{k-1},m_k,n_{k-1},n_k)\in\ZZ^4 : \\
|m_{k-1}|,|m_k|,|n_{k-1}|,|n_k|\ll X^{1/k},\;(m_{k-1},m_k)\neq (n_{k-1},n_k), \\
(m_{k-1}-\alpha_{k-1})^2 + (m_k-\alpha_k)^2
=(n_{k-1}-\alpha_{k-1})^2 + (n_k-\alpha_k)^2 \} \\
\sim X^{(k-2)/k}\times \tilde c_\vecalf X^{2/k} \log X ,
\end{multline*}
as $X\rightarrow\infty$, for some constant $\tilde c_\vecalf>0$
(compare Section 9 in \cite{Marklof00}). We conclude that,
for $X$ large enough,
\begin{multline*}
\frac1X\# \{ (\vecm,\vecn) \in\ZZ^k\times\ZZ^k: \; \vecm\neq\vecn,\; 
\|\vecm-\vecalf\|^k \leq X, \|\vecn-\vecalf\|^k \leq X, \; \\
 \|\vecm-\vecalf\|^k=\|\vecn-\vecalf\|^k \}
\geq c_\vecalf  \log X ,
\end{multline*}
for some constant $c_\vecalf>0$.

\subsection{\label{bubu2}}
By a similar argument, one has for $\vecalf\in\QQ^k$
\begin{multline*}
\frac1X\# \{ (\vecm,\vecn) \in\ZZ^k\times\ZZ^k: \; \vecm\neq\vecn,\; 
\|\vecm-\vecalf\|^k \leq X, \|\vecn-\vecalf\|^k \leq X, \; \\
 \|\vecm-\vecalf\|^k=\|\vecn-\vecalf\|^k \}
\sim c_\vecalf X^{(k-2)/k}  
\end{multline*}
for $X\rightarrow\infty$.
This can be derived, e.g.~in the case $\vecalf=\vecnull$, from the asymptotics
$$
\int_0^1
\Theta_f(u+\i \frac{1}{\lambda},0;\vecnull)
\overline{\Theta_g(u+\i \frac{1}{\lambda},0;\vecnull)} \; du
\sim b \lambda^{k/2-1},
$$
compare, e.g., Theorem 6.1 \cite{Marklof99}.

\subsection{Proof of Theorem \ref{sharp}}
Let ${\mathcal B}$ be a countable dense
set of badly approximable $(k-2)$-tuples.
Enumerate the quadratic forms $\|\vecx-\vecalf_j\|^2$ with
$\vecalf_j\in{\mathcal B}\times\QQ^2$ as $P_1,P_2,P_3,\ldots$. Because of the
bound derived in \ref{bubu}, given any $X>1$, there exists
an $M_j>X$ such that 
$$
R_2^{\vecalf_j}[0,0](M_j) \geq \frac{\log M_j}{\log\log\log M_j} .
$$
We find a small $\epsilon_j=\epsilon_j(M_j)>0$ such that
$$
R_2^{\vecalf}[-a,a](M_j)\geq R_2^{\vecalf_j}[0,0](M_j) 
$$
for all $\vecalf\in B_j$, where $B_j$ is the open set of all
$\vecalf$ with $\|\vecalf-\vecalf_j\|<\epsilon_j$.
Individually, the sets $B_j$ shrink
to a point as $X\rightarrow\infty$,
but the union
$$
\bigcup_{j: M_j\geq X} B_j
$$
is open and dense in $\TT^k$. Therefore
$$
B=\bigcap_{X=1}^{\infty} \bigcup_{j: M_j\geq X} B_j
$$
is of second Baire category.
So if $\vecalf\in B$, then, given any $X$, 
there exists some $M\geq X$, such that
$$
R_2^{\vecalf}[-a,a](M)\geq \frac{\log M}{\log\log\log M}.
$$

Note that the proof remains valid if $\log\log\log$ is replaced by any
slowly increasing positive function $\nu \leq \log\log\log$
with $\nu(X)\rightarrow\infty$ ($X\rightarrow\infty$).

Property (iii) follows from Theorem \ref{pairthm} by the same 
string of arguments used in Section 9.3 \cite{Marklof00}.
\qed

\subsection{Proof of Theorem \ref{except}}
Follows from the relation in Section \ref{bubu2}. The proof is otherwise
identical to the proof of Theorem 1.13 \cite{Marklof00}.
\qed


\begin{thebibliography}{99}

\bibitem{Berry77}
M.V.\,Berry and M.\,Tabor, 
Level clustering in the regular spectrum,
{\em Proc. Roy. Soc.} A {\bf 356} (1977) 375-394.

\bibitem{Bleher99}
P.M.\,Bleher,
Trace formula for quantum integrable systems, 
lattice point problem, and small divisors,
in: D.\,Hejhal et al.~(eds.),
{\em Emerging Applications of Number Theory},
IMA Volumes in Mathematics and its Applications,
Vol.~109 (Springer, New York, 1999) pp.~1-38. 

\bibitem{Cheng91}
Z.\,Cheng and J.L.\,Lebowitz, Statistics of energy levels in integrable
quantum systems,
{\em Phys. Rev.} A {\bf 44} (1991) 3399-3402.

\bibitem{Cheng94}
Z.\,Cheng, J.L.\,Lebowitz and P.\,Major, On the number of lattice
points between two enlarged and randomly shifted copies of an oval, 
{\em Probab. Theory Related Fields} {\bf 100} (1994) 253-268. 

\bibitem{Dani93}
S.G.\,Dani and G.A.\,Margulis, 
Limit distributions of orbits of unipotent flows and values of 
quadratic forms, {\em Adv. Soviet Math.} 16, Part 1, 
(Amer. Math. Soc., Providence, RI, 1993)  91-137.

\bibitem{Eskin98b}
A.\,Eskin, G.\,Margulis and S.\,Mozes, 
Quadratic forms of signature (2,2) 
and eigenvalue spacings on flat 2-tori, preprint.

\bibitem{Lion80}
G.\,Lion and M.\,Vergne,
{\em The Weil Representation, Maslov Index and Theta
Series}, Progr. in Math. 6 (Birkh\"auser, 1980).

\bibitem{Marklof98}
J.\,Marklof, Spectral form factors of rectangle billiards,
{\em Comm. Math. Phys.} {\bf 199} (1998) 169-202.

\bibitem{Marklof99}
J.\,Marklof, Limit theorems for theta sums,
{\em Duke Math. J.} {\bf 97} (1999) 127-153.

\bibitem{Marklof00}
J.\,Marklof, Pair correlation densities of inhomogeneous quadratic forms,
{\em Ann. of Math.}, to appear.

\bibitem{Marklof00c}
J.\,Marklof,
The Berry-Tabor conjecture,
{\em Proceedings of the 3rd European Congress of Mathematics},
Barcelona 2000, Progress in Mathematics Vol.~202
(Birkh\"auser, Basel, 2001) 421-427.

\bibitem{Marklof00d}
J.\,Marklof,
Level spacing statistics and integrable dynamics,
{\em X\nolinebreak I\nolinebreak I\nolinebreak Ith 
International Congress on Mathematical Physics},
London 2000 (International Press, Boston, 2001) 359-363.

\bibitem{Sarnak97}
P.\,Sarnak, Values at integers of binary quadratic forms, 
{\em Harmonic Analysis and Number Theory} (Montreal, PQ, 1996), 
181-203, CMS Conf. Proc. {\bf 21}, Amer. Math. Soc., Providence, RI, 1997. 

\bibitem{Sarnak97I}
P.\,Sarnak,
Quantum chaos, symmetry and zeta functions.
Lecture I: Quantum chaos,
{\em Curr. Dev. Math.} (1997) 84-101.

\bibitem{Schmidt80}
W.M.\,Schmidt,
{\em Diophantine Approximation},
Lect. Notes Math. {\bf 785} (Springer, 1980).

\bibitem{Shah96}
N.A.\,Shah, Limit distributions of expanding 
translates of certain orbits on homogeneous spaces,
{\em Proc. Indian Acad. Sci., Math. Sci.} {\bf 106} (1996) 105-125.

\bibitem{VanderKam96}
J.M.\,VanderKam, Values at integers of homogeneous polynomials,
{\em Duke Math. J.} {\bf 97} (1999) 379-412.

\bibitem{VanderKam97}
J.M.\,VanderKam, Pair correlation of four-dimensional flat tori,
{\em Duke Math. J.} {\bf 97} (1999) 413-438.

\bibitem{VanderKam00}
J.M.\,VanderKam, Correlations of eigenvalues on multi-dimensional flat tori,
{\em Comm. Math. Phys.} {\bf 210} (2000) 203-223.

\end{thebibliography}
\end{document}